\documentclass[12pt]{article}
\usepackage{amsfonts,amsmath}
\usepackage{epsfig,graphics}

\renewcommand{\appendix}{%
\renewcommand{\section}{%
\newpage
\thispagestyle{plain}%
\secdef\Appendix\sAppendix}%
\setcounter{section}{0}%
\renewcommand{\thesection}{\Alph{section}}%
}
\newcommand{\Appendix}[2][?]{%
\refstepcounter{section}%
\addcontentsline{toc}{Addendum}%
{\protect\numberline{\appendixname~\thesection}#1}%
{\flushleft\LARGE\bfseries\appendixname\ \thesection\par
\centering#2\par}%
\sectionmark{#1}\vspace{\baselineskip}}

\newcommand{\sAppendix}[1]{%
{\flushright\large\bfseries\appendixname\par
\centering#1\par}%
\vspace{\baselineskip}}

\def\be{\begin{equation}}

\def\bea{\begin{eqnarray}}

\def\eea{\end{eqnarray}}

\normalbaselines \oddsidemargin 0.5cm \evensidemargin 0.5cm
\begin{document}

\pagestyle{empty}

\rightline{RR.047.0607}

\vskip 1cm

\begin{center}

{\Huge {\textbf{Higher Dimensional Multiparameter Unitary and
Nonunitary Braid Matrices: Even Dimensions}}}

\vspace{4mm}

{\bf \large B. Abdesselam$^{a,}$\footnote{Email:
boucif@cpht.polytechnique.fr and  boucif@yahoo.fr}, A.
Chakrabarti$^{b,}$\footnote{Email:
chakra@cpht.polytechnique.fr},\\
V.K. Dobrev$^{c,d,}$\footnote{Email: dobrev@inrne.bas.bg} and S.G.
Mihov$^{c,}$\footnote{Email: smikhov@inrne.bas.bg}}

\vspace{2mm}

  \emph{$^a$ Laboratoire de Physique Quantique de la
Mati\`ere et de Mod\'elisations Math\'ematiques, Centre
Universitaire de Mascara, 29000-Mascara, Alg\'erie\\
and \\
Laboratoire de Physique Th\'eorique, Universit\'e d'Oran
Es-S\'enia, 31100-Oran, Alg\'erie}
  \\
  \vspace{2mm}
  \emph{$^b$ Centre de Physique Th{\'e}orique, Ecole Polytechnique, 91128 Palaiseau Cedex, France.}
  \\
  \vspace{2mm}
  \emph{$^c$ Institute of Nuclear Research and Nuclear Energy
Bulgarian Academy of Sciences 72 Tsarigradsko Chaussee, 1784
Sofia, Bulgaria}\\
 \vspace{2mm}
 \emph{$^d$ Abdus Salam International Center for Theoretical Physics
Strada Costiera 11, 34100 Trieste, Italy}

\end{center}

\begin{abstract}
{\small \noindent A class of $\left(2n\right)^2\times\left(
2n\right)^2$ multiparameter braid matrices are presented for all
$n$ $\left(n\geq 1\right)$. Apart from the spectral parameter
$\theta$, they depend on $2n^2$ free parameters $m_{ij}^{(\pm)}$,
$i,j=1,\ldots,n$. For real parameters the matrices
$R\left(\theta\right)$ are nonunitary. For purely imaginary
parameters they became unitary. Thus a unification is achieved
with odd dimensional multiparameter solutions presented before.}
\end{abstract}

\pagestyle{plain} \setcounter{page}{1}



\section{Introduction}
\setcounter{equation}{0}

Higher dimensional unitary braid matrices have been studied in two
recent papers \cite{R1,R2}. Their simultaneous relevance to
topological and quantum entanglements (as discussed, for example,
in Ref. \cite{R3}) was a major motivation. In Ref. \cite{R2} quite
different classes were presented for odd and even dimensional
matrices. There, the even dimensional $\left(2n\right)^2\times
\left(2n\right)^2$ braid matrices have no free parameter (apart
from the spectral parameter $\theta$ after Baxterization) where as
the $\left(2n+1\right)^2\times \left(2n+1\right)^2$ matrices have
$2n\left(n+2\right)$ free parameters
$\left(m_{ij}^{(\pm)}\right)$. Here we unify the two cases by
presenting multiparameter solutions for even dimensions. We obtain
first the general case for this class and then show how to
implement unitarity.

\section{Constructions (Even dimensions)}
\setcounter{equation}{0}

The braid equation is, in standard notations, in presence of a
spectral parameter $\theta$,
\begin{equation}
\widehat{R}_{12}\left(\theta\right)\widehat{R}_{23}\left(\theta+\theta'\right)\widehat{R}_{12}\left(\theta'\right)=
\widehat{R}_{23}\left(\theta'\right)\widehat{R}_{12}\left(\theta+\theta'\right)\widehat{R}_{23}\left(\theta\right),
\end{equation}
where $\widehat{R}_{12}=\widehat{R}\otimes I$ and
$\widehat{R}_{23}=I\otimes \widehat{R}$. We present below a simple
class of multiparameter solutions for $\left(2n\right)^2\times
\left(2n\right)^2$ $\left(n\geq 1\right)$ braid matrices
$\hat{R}\left(\theta\right)$. They are analogous to the odd
dimensional solutions presented before \cite{R4}. Unitarity
constraints can be implemented as in sec. 5 of Ref. \cite{R2}.
Thus, for this class, one obtains a unified approach for
multiparameter odd and even dimensional solutions.

Define the projectors
\begin{equation}
P_{ij}^{(\epsilon)}=\frac 12\left\{\left(ii\right)\otimes
\left(jj\right)+\left(\bar{i}\bar{i}\right)\otimes
\left(\bar{j}\bar{j}\right)+\epsilon\left[\left(i\bar{i}\right)\otimes
\left(j\bar{j}\right)+\left(\bar{i}i\right)\otimes
\left(\bar{j}j\right)\right]\right\},
\end{equation}
where $i,j\in\left\{1,\ldots,n\right\}$, $\bar{i}=2n+1-i$,
$\bar{j}=2n+1-j$ and $\epsilon=\pm$. They provide a complete basis
satisfying
\begin{equation}
P_{ij}^{(\epsilon)}P_{kl}^{(\epsilon')}=\delta_{ik}\delta_{jl}\delta_{\epsilon\epsilon'}P_{ij}^{(\epsilon)},\qquad
\sum_{\epsilon=\pm}\sum_{i,j=1}^nP_{ij}^{(\epsilon)}=I_{(2n)^2\times
(2n)^2}.
\end{equation}
Anticipating the basic constraints essential for odd dimension
\cite{R4} we directly postulate the form
\begin{equation}
\hat{R}\left(\theta\right)=\sum_{\epsilon=\pm}\sum_{i,j=1}^ne^{m_{ij}^{(\epsilon)}\theta}\left(P_{ij}^{(\epsilon)}+P_{i\bar{j}}^{(\epsilon)}\right).
\end{equation}
The proof that it satisfies the braid equation (2.1) proceeds in
close analogy to the equations from (A9) to (A17) of Ref.
\cite{R4}. Here we have directly implemented the constraint
\begin{equation}
m_{ij}^{(\epsilon)}=m_{i\bar{j}}^{(\epsilon)}.
\end{equation}
It is instructive to study explicitly the simplest cases.

\begin{description}
\item[Case 1:] $N=2$ ($n=1$) (Here $i=1$, $\bar{i}=2$ and
similarly for $j$).
\begin{equation}
\begin{pmatrix}
  a_+ & 0 & 0 & a_- \\
  0 & a_+ & a_- & 0 \\
  0 & a_- & a_+ & 0 \\
  a_- & 0 & 0 & a_+ \\
\end{pmatrix},
\end{equation}
with
\begin{equation}
a_{\pm}=\frac 12\left(e^{m_{11}^{(+)}\theta}\pm
e^{m_{11}^{(-)}\theta}\right).
\end{equation}

\item[Case 2:] $N=4$ ($n=2$) (Here $i=1,2$, $\bar{i}=3,4$). In
terms of $4\times 4$ blocks ($D_{ij}$, $A_{ij}$ on the diag. and
anti-diag. respectively)
\begin{equation}
\begin{pmatrix}
  D_{11} & 0 & 0 & A_{1\bar{1}} \\
  0 & D_{22} & A_{2\bar{2}} & 0 \\
  0 & A_{\bar{2}2} & D_{\bar{2}\bar{2}} & 0 \\
  A_{\bar{1}1} & 0 & 0 & D_{\bar{1}\bar{1}} \\
\end{pmatrix},
\end{equation}
with
\begin{eqnarray}
&&D_{11}=D_{\bar{1}\bar{1}}=\left(\begin{array}{cccc}
  a_+ & 0 & 0 & 0 \\
  0 & b_+ & 0 & 0 \\
  0 & 0 & b_+ & 0 \\
  0 & 0 & 0 & a_+ \\
\end{array}\right),\qquad
D_{22}=D_{\bar{2}\bar{2}}=\left(\begin{array}{cccc}
  c_+ & 0 & 0 & 0 \\
  0 & d_+ & 0 & 0 \\
  0 & 0 & d_+ & 0 \\
  0 & 0 & 0 & c_+ \\
\end{array}\right),\nonumber\\
&&A_{1\bar{1}}=A_{\bar{1}1}=\left(\begin{array}{cccc}
  0 & 0 & 0 & a_- \\
  0 & 0 & b_- & 0 \\
  0 & b_- & 0 & 0 \\
  a_- & 0 & 0 & 0 \\
\end{array}\right),\qquad
A_{2\bar{2}}=A_{\bar{2}2}=\left(\begin{array}{cccc}
  0 & 0 & 0 & c_- \\
  0 & 0 & d_- & 0 \\
  0 & d_- & 0 & 0 \\
  c_- & 0 & 0 & 0 \\
\end{array}\right)
\end{eqnarray}
and
\begin{eqnarray}
&&a_{\pm}=\frac 12\left(e^{m_{11}^{(+)}\theta}\pm
e^{m_{11}^{(-)}\theta}\right),\qquad b_{\pm}=\frac
12\left(e^{m_{12}^{(+)}\theta}\pm e^{m_{12}^{(-)}\theta}\right),\nonumber\\
&&c_{\pm}=\frac 12\left(e^{m_{21}^{(+)}\theta}\pm
e^{m_{21}^{(-)}\theta}\right),\qquad d_{\pm}=\frac
12\left(e^{m_{22}^{(+)}\theta}\pm
e^{m_{22}^{(-)}\theta}\right)\end{eqnarray}

\end{description}
We have verified, using a program, the braid equation (2.1) by
inserting (2.4) for $N=2,4,6,8$. These provide direct checks for
the argument indicated below (2.4). As compared to $\frac
12\left(N+3\right)\left(N-1\right)$ free parameters for odd
\cite{R4}, here for even $N$ we obtain $\frac 12N^2$ free
parameters $m_{ij}^{(\pm)}$.

Let us just note that the odd dimensional solutions of Ref.
\cite{R4} and the even dimensional solutions presented here can be
regrouped in a single expression given by
\begin{equation}
{\hat R}\left(\theta\right)=\frac 12\sum_{\epsilon=\pm}
\sum_{i,j=1}^N
e^{m_{ij}^{(\epsilon)}\theta}\left[\left(ii\right)\otimes
\left(jj\right)+\epsilon\left(i{\bar i}\right)\otimes \left(j{\bar
j}\right)\right],
\end{equation}
where
\begin{eqnarray}
&&m_{ij}^{(\epsilon)}=m_{{\bar i}j}^{(\epsilon)}=m_{i{\bar
j}}^{(\epsilon)}=m_{{\bar i}{\bar j}}^{(\epsilon)},\qquad \qquad
i,j=1,\ldots,N\,\hbox{and}\, \epsilon=\pm 1,\nonumber \\
&& \overline{n+1}=n+1\,\,\hbox{and}\,\,
m_{n+1,n+1}^{(\epsilon)}=0\,(\forall\epsilon)\qquad\qquad \hbox{If
$N$ is odd, i.e. $N=2n+1$}.
\end{eqnarray}

\section{Unitarity}
\setcounter{equation}{0}

For all parameter real, $\hat{R}\left(\theta\right)$ is real but
not unitary. Exactly as for $N$ odd, making each exponent purely
imaginary, namely
$\exp\left({m_{ij}^{(\pm)}\theta}\right)\longrightarrow
\exp\left({{\mathbf i}m_{ij}^{(\pm)}\theta}\right)$, where on the
right $m_{ij}^{(\pm)}\theta$ is now real with a coefficient
${\mathbf i}$ $\left({\mathbf i}^2=-1\right)$, one obtains
unitarity. Now, due to the symmetry of the projectors
\begin{equation}
\hat{R}\left(\theta\right)^{+}=\hat{R}\left(-\theta\right)=\hat{R}\left(\theta\right)^{-1},\qquad
\hat{R}\left(\theta\right)^{+}\hat{R}\left(\theta\right)=I_{(2n)^2\times
(2n)^2}.
\end{equation}
In general, one can demonstrate that our multiparameter odd and
even dimensional solutions one has a simple factorization
\begin{equation}
\hat{R}\left(\theta_1\pm\theta_2\right)=\hat{R}\left(\theta_1\right)\hat{R}\left(\pm\theta_2\right)=
\hat{R}\left(\theta_1\right)\hat{R}\left(\theta_2\right)^{\pm 1}.
\end{equation}
This evidently, holds for real or imaginary parameters, i.e. for
nonunitary and unitary solutions. Correspondingly, the $R{\mathbf
{TT}}$ relations can be expressed as follows:
\begin{equation}
\left(\hat{R}\left(\theta\right)\left({\mathbf
T}\left(\theta\right)\otimes I\right)\right)\left(\left(I\otimes
{\mathbf
T}\left(\theta'\right)\right)\hat{R}\left(\theta'\right)\right)
=\left(\hat{R}\left(\theta'\right)\left({\mathbf
T}\left(\theta'\right)\otimes I\right)\right)\left(\left(I\otimes
{\mathbf
T}\left(\theta\right)\right)\hat{R}\left(\theta\right)\right).
\end{equation}

For comparison one may note that the real unitary braid matrix for
all $N$ ($=2n$) presented in \cite{R2} can be written as
\begin{equation}
\hat{R}\left(z\right)=\displaystyle\left(\frac{1-{\mathbf
i}z}{1+{\mathbf i}z}\right)^{1/2}P_+ +\left(\frac{1+{\mathbf
i}z}{1-{\mathbf i}z}\right)^{1/2}P_-,
\end{equation}
where
\begin{equation}
z=\tanh\left(\theta\right),\qquad \left(\frac{1\mp {\mathbf
i}z}{1\pm {\mathbf i} z}\right)^{1/2}\equiv e^{\pm {\mathbf
i}\varphi},
\end{equation}
say, giving phases for the coefficients and
\begin{equation}
P_{\pm}=\frac 12\left(I\otimes I\pm {\mathbf i}K\otimes J\right).
\end{equation}
($K$, $J$ being given by (2.2) of Ref. \cite{R2}). $P_{\pm}$ can
be expressed as sums of the projectors of the type
\begin{equation}
Q_{ij}^{(\epsilon)}=\frac 12\left\{\left(ii\right)\otimes
\left(jj\right)+\left(\bar{i}\bar{i}\right)\otimes
\left(\bar{j}\bar{j}\right)+\epsilon {\mathbf
i}(-1)^{\bar{j}}\left[\left(i\bar{i}\right)\otimes
\left(j\bar{j}\right)-\left(\bar{i}i\right)\otimes
\left(\bar{j}j\right)\right]\right\}
\end{equation}
defining analogously $Q_{i\bar{j}}^{(\epsilon)}$ (with
$j\longrightarrow \bar{j}$, $\bar{j}\longrightarrow j$ in
$Q_{i\bar{j}}^{(\epsilon)}$). The imaginary factor ${\mathbf i}$
in $Q_{ij}^{(\epsilon)}$, $Q_{i\bar{j}}^{(\epsilon)}$ and the
phases cancel giving a real $\hat{R}\left(z\right)$,
\begin{equation}
\hat{R}\left(z\right)=I\otimes I+ zK\otimes J.
\end{equation}
Due to the summing up of $Q_{ij}^{(\epsilon)}$ into $P_{\pm}$ the
effective number of projectors does not increase with $N$, nor the
number of parameters. Here we have presented a class of solutions
where the number of free parameters increase as $N^2$. For this
case, one can prove that
\begin{equation}
\hat{R}\left(z_1\right)\hat{R}\left(z_2\right)=\hat{R}\left(z_3\right),
\end{equation}
where
\begin{equation}
z_3=\frac{z_1+z_2}{1-z_1z_2}\neq
\tanh\left(\theta_1+\theta_2\right),\qquad \left(z_1z_2\neq
1\right).
\end{equation}

\section{$\theta$-Expansion}
\setcounter{equation}{0}

In section 5 of Ref. \cite{R4} exponentiation and
$\theta$-expansion of $\hat{R}\left(\theta\right)$ was presented
for odd dimension. We present below a brief analogous treatment
for even dimensions. One have
\begin{equation}
e^{m_{ij}^{(\epsilon)}\theta}\left(P_{ij}^{(\epsilon)}+P_{i\bar{j}}^{(\epsilon)}\right)=
\left(P_{ij}^{(\epsilon)}+P_{i\bar{j}}^{(\epsilon)}\right)+\sum_{k=1}^{\infty}\frac
1{k!}\left(m_{ij}^{(\epsilon)}\left(P_{ij}^{(\epsilon)}+P_{i\bar{j}}^{(\epsilon)}\right)\right)^k\theta^k.
\end{equation}
Defining
\begin{equation}
X=\sum_{\epsilon=\pm}\sum_{i,j=1}^nm_{ij}^{(\epsilon)}\left(P_{ij}^{(\epsilon)}+P_{i\bar{j}}^{(\epsilon)}\right)\Longrightarrow
X^n=\sum_{\epsilon=\pm}\sum_{i,j=1}^n\left(m_{ij}^{(\epsilon)}\right)^n\left(P_{ij}^{(\epsilon)}+P_{i\bar{j}}^{(\epsilon)}\right)
\end{equation}
due to the orthogonality of the projectors
$\left(P_{ij}^{(\epsilon)}+P_{i\bar{j}}^{(\epsilon)}\right)$ for
different sets of indices. Now from (2.4), due to the completeness
(2.3),
\begin{equation}
\hat{R}\left(\theta\right)=\sum_{\epsilon=\pm}\sum_{i,j=1}^ne^{m_{ij}^{(\epsilon)}\theta}\left(P_{ij}^{(\epsilon)}+P_{i\bar{j}}^{(\epsilon)}\right)=
I+\sum_{k=1}^{\infty}\frac 1{k!}X^k\theta^k=e^{\theta X}.
\end{equation}
Hence the braid equation (2.1) reduces to
\begin{equation}
e^{\theta X_{12}}e^{\left(\theta+\theta'\right) X_{23}}e^{\theta'
X_{12}}=e^{\theta' X_{23}}e^{\left(\theta+\theta'\right)
X_{12}}e^{\theta X_{23}},
\end{equation}
where $X_{12}=X\otimes I$ and $X_{23}=I\otimes X$. Expanding both
sides and comparing coefficients of $\theta^{a}\left(\theta+
\theta'\right)^b\theta'^c$ one obtains a sequence of relations
involving $X_{12}$, $X_{23}$. Some have been pointed out in
section 5 of Ref. \cite{R4}. There would be parallel features
here.

After implementing unitarity as in section 3 one can define (with
${\mathbf i}$ as above (3.1))
\begin{eqnarray}
X={\mathbf
i}\sum_{i,j,\epsilon}m_{ij}^{(\epsilon)}\left(P_{ij}^{(\epsilon)}+P_{i\bar{j}}^{(\epsilon)}\right).
\end{eqnarray}
One then proceeds as above.

We have started with multiparameter case. For
\begin{equation}
\hat{R}\left(z\right)=I\otimes I+ zK\otimes
J=\left(\frac{1-{\mathbf i}z}{1+{\mathbf i}z}\right)^{1/2}P_+
+\left(\frac{1+{\mathbf i}z}{1-{\mathbf i}z}\right)^{1/2}P_-\equiv
e^{{\mathbf i}\varphi}P_++ e^{-{\mathbf i}\varphi}P_-.
\end{equation}
(see the discussion following (3.3))
\begin{eqnarray}
X={\mathbf i}\left(P_{+}-P_{-}\right)=-K\otimes J
\end{eqnarray}
and
\begin{equation}
\hat{R}\left(z\right)=\hat{R}\left(\varphi\right)=e^{\varphi X},
\end{equation}
where $e^{{\mathbf i}\varphi}=\left(\frac{1-{\mathbf
i}\tanh\theta}{1+{\mathbf i}\tanh\theta}\right)^{1/2}$. By
inserting this $X$ in (4.4), one can develop in $\varphi$ as
explained before.

\section{Discussion}
\setcounter{equation}{0}

For the multiparameter solutions presented in sections 2 and 3 one
can study $\hat{R}\verb"T"\verb"T"$ relations, transfer matrices,
Hamiltonians and factorizable $S$-matrices in a closely analogous
fashion to that for odd dimensions \cite{R5}. They are beyond the
scope of this paper, limited essentially to construction of
multiparameter $\left(2n\right)^2\times\left(2n\right)^2$ braid
matrices (nonunitary and unitary).

Beyond the unification presented there is a basic difference
between odd and even dimensional cases. For the
$\left(2n+1\right)^2\times \left(2n+1\right)^2$ braid matrices
with a basis of our "nested sequence" of projectors our
multiparameter solutions are the most general ones. The presence
of the central element 1 in $\hat{R}\left(\theta\right)$ imposes
the simple exponential solutions for the coefficients of all other
projectors. This has been emphasized in appendix A of Ref.
\cite{R4}. ("Solving the braid equation"). But for even dimension
this not the case. The class of solutions presented here is only
one possibility. Already for the $4\times 4$ case the intensively
studied 6- and 8-vertex solutions can be canonically expressed on
our basis (sections 6 and 7 of Ref. \cite{R6}). The
multidimensional generalization of the 6-vertex matrix presented
in Ref. \cite{R7} (citing previous sources) remains restricted to
a single parameter $\gamma$. Are there authentic multiparameter
generalizations of 6- and 8-vertex solutions to $\left(2n
\right)^2\times\left(2n \right)^2$ matrices for $n>1$? We intend
to explore such possibilities elsewhere.

We point out moreover that a pure imaginary spectral parameter
$\left(\theta\longrightarrow {\mathbf i}\theta\right)$ renders the
6-vertex and 8-vertex braid matrices unitary. This particularly
evident form the respective canonical forms ((6.5) for 6-vertex
and (7.2) with (7.6), (7.7) for 8-vertex of Ref. \cite{R6}) where
the normalization factors are also suitably adapted. The
coefficient of each real symmetric projector is evidently inverted
under conjugation $\left({\mathbf i}\theta\longrightarrow
-{\mathbf i}\theta\right)$. Hence
\begin{equation}
\hat{R}^{+}\left(\theta\right)\hat{R}\left(\theta\right)=\hat{R}\left(-\theta\right)\hat{R}\left(\theta\right)=I.
\end{equation}
Now one no longer has statistical models with real, non-negative
Boltzmann weights. But the unitary matrices become relevant
concerning entanglement. Such parametrizations of entangled states
will be studied in a following paper.

The unitary $\left(2n \right)^2\times\left(2n \right)^2$ braid
matrices generate entangled quantum states with one difference as
compared to the odd dimensional case. In the last section of Ref.
\cite{R2} it was pointed out that the product of pure states
$\left|0\right\rangle\left|0\right\rangle$ conserved its status
under action of $\left(2n+1\right)^2\times\left(2n+1\right)^2$
unitary matrix. For the present case there is no such exceptional
state.

As already pointed out for odd dimensions (see section 5, Ref.
\cite{R2}), even dimensional, unitary, multiparameter braid
matrices are also periodic or quasiperiodic in $\theta$
accordingly as the $m$'s are mutually commensurate or not.

\vskip 0.5cm

\noindent{\bf Acknowledgments:} {\em One of us (BA) wants to thank
Pierre Collet and Paul Sorba for precious help. The work of VKD
and SGM was supported in part by the Bulgarian National Council
for scientific Research, grant F-1205/02 and the European RTN
'Force-universe', contract MRTN-CT-2004-005104, and by the
Alexander Von Humboldt Foundation in framework of the
Clausthal-Leipzig-Sofia cooperation. }

\end{document}